\index{\right] }
\documentclass{amsart}
\usepackage{amsthm}
\usepackage{amsfonts}
\usepackage{amssymb}
\usepackage{amsmath}

\begin{document}
\small \setcounter {page}{1}

\title{\textbf{STATISTICAL CONVERGENCE OF ORDER $\alpha$ IN PROBABILITY}}


\author{Pratulananda Das$^{(1)}$, Sanjoy Ghosal$^{(2)}$, Sumit Som$^{(3)*}$  }

\thanks{$^{*}$ This work is funded by UGC Research, HRDG, India. $^{(1)}$Department of Mathematics, Jadavpur University, Kolkata-700032, West Bengal, India.  E-mail:  pratulananda@yahoo.co.in  (P. Das).
$^{(2)}$Department of Mathematics, Kalyani Government Engineering College, Kalyani, Nadia-741235, West Bengal, India. E-mail:  sanjoykrghosal@yahoo.co.in  (S. Ghosal).
$^{(3)}$UGC Fellow, Department of Mathematics, Jadavpur University, Kolkata-700032, West Bengal, India. E-mail:  somkakdwip@gmail.com  (S. Som)}
\date{}

\setcounter {page}{1}

\noindent{\maketitle{\bf Abstract:}} In this paper ideas of different types of convergence of a sequence of random variables in probability, namely, statistical convergence of order $\alpha$ in probability, strong $p$-Ces$\grave{\mbox{a}}$ro summability of order $\alpha$ in probability, lacunary statistical convergence or $S_{\theta}$-convergence of order $\alpha$ in probability, ${N_{\theta}}$-convergence of order $\alpha$ in probability have been introduced and their certain basic properties have been studied.\\



\section{\textbf{Introduction and Background}}
\baselineskip .82cm \small

The idea of convergence of a real sequence has been extended to statistical convergence by Fast \cite{6} and Steinhaus \cite{24} as follows: If $\mathbb{N}$ denotes the set of natural numbers and $K\subset \mathbb{N},$ then $K(m,n)$ denotes the cardinality of the set $K\cap [m,n].$ The upper and lower natural density of the subset $K$ is defined by $$\bar{d}(K)=\displaystyle{\lim_{n\rightarrow \infty}}\sup \frac{K(1,n)}{n} ~ ~ \mbox{and} ~ ~ \underline{d}(K)=\displaystyle{\lim_{n\rightarrow \infty}}\inf \frac{K(1,n)}{n}.$$
If $\bar{d}(K)=\underline{d}(K),$ then we say that the natural density of $K$ exists, and it is denoted simply by $$d(K)=\displaystyle{\lim_{n\rightarrow \infty}} \frac{K(1,n)}{n}.$$
A sequence $\{x_{n}\}_{n\in \mathbb{N}}$  of real numbers is said to be statistically convergent to a real number $x$  if for each $\varepsilon>0,$ the set $K=\{n\in \mathbb{N}: |x_n-x|\geq \varepsilon\}$ has natural density zero and we write $x_{n}\xrightarrow{S} x.$  Statistical convergence has turned out to be one of the most active areas of research in summability theory after the work of Fridy \cite{7} and \v{S}al\'{a}t \cite{21}. Over the years a lot of work have been done to generalize this notion of statistical convergence and to introduce new summability methods related to it. Some of the most important concepts introduced are : lacunary statistical convergence by Fridy \& Orhan \cite{9} (for more results on this convergence see the paper of Li \cite{16}), $\mathcal{I}$-convergence by Kostyrko et al. \cite{13}, (see \cite{1,D1,D2,4,22} for recent advances and more references on this convergence), statistical convergence of order $\alpha$ by Bhunia et. al. \cite{2} (statistical convergence of order $\alpha$ was also independently introduced by Colak \cite{3}, more investigations in this direction and more applications can be found in \cite{C2}), lacunary statistical convergence of order $\alpha$ by Seng$\ddot{u}$l \&  Et. M \cite{E1},
  pointwise and uniform statistical convergence of order $\alpha$ for sequences of functions by Cinar et al. \cite{E2}, $\lambda$-statistical convergence of order $\alpha$ of sequences of function by Et. M et al. \cite{E3}, $\mathcal{I}$-statistical and $\mathcal{I}$-lacunary statistical convergence of order $\alpha$ by Savas \& Das \cite{23}, open covers and selection principles by Das \cite{pd1,pd2}. The notion of statistical convergence has applications in different fields of mathematics : in number theory by Erd$\ddot{\mbox{o}}$s \& Tenenbaum \cite{5}, in statistics and probability theory by Fridy \& Khan \cite{8} and Ghosal \cite{11, 12, G3}, in approximation theory by Gadjiev \& Orhan \cite{10}, in hopfield neural network by Martinez et al. \cite{17}, in optimization by Pehlivan \&  Mamedov \cite{19}.

In particular in probability theory, if for each positive integer $n$, a random variable $X_n$ is defined on a given event space $S$ (same for each $n$) with respect to a given class of events $\triangle$ and a probability function $P:\triangle \rightarrow \mathbb{R}$ (where $\mathbb{R}$ denotes the set of real numbers) then we say that $X_1,X_2,X_3,...,X_n,...$ is a sequence of random variables and as in analysis we denote this sequence by $\{X_{n}\}_{n\in \mathbb{N}}.$

From the practical point of view the discussion of a random variable $X$ is highly significant if it is known that there exists a real constant $c$ for which $P(|X-c|<\epsilon)\simeq 1,$ where $\epsilon>0$ is sufficiently small, that is, it is nearly certain that values of $X$ lie in a very small neighbourhood of $c.$

For a sequence of random variables $\{X_{n}\}_{n\in \mathbb{N}},$ each $X_{n}$ may not have the above property but it may happen that the aforementioned property (with respect to a real constant $c$) becomes more and more distinguishable as $n$ gradually increases and the question of existence of such a real constant $c$ can be answered by a concept of convergence in probability of the sequence $\{X_{n}\}_{n\in \mathbb{N}}.$

        In this short paper we shall limit our discussion to four types of convergence of a sequence of random variables, namely,\\
      (i)  statistical convergence of order $\alpha$ in probability,\\
      (ii) strong $p$-Ces$\grave{\mbox{a}}$ro summability of order $\alpha$ in probability,\\
      (iii) lacunary statistical convergence or $S_{\theta}$-convergence of order $\alpha$ in probability,\\
      (iv) ${N_{\theta}}$-convergence of order $\alpha$ in probability.\\

Our main aim in this paper is to establish some important theorems related to the modes of convergence (i) to (iv), which effectively extend and improve all the existing results in this direction \cite{2,3,9,11,16,20}. We also intend to establish the relation between these four summability notions. It is important to note that the methods of proofs and in particular the examples are not analogous to the real case. \\

\section{\textbf{Statistical Convergence of Order $\alpha$ in Probability}}

We first recall the definition of statistical convergence of order $\alpha$ of a sequence of real numbers from (\cite{2,3}) as follows:\\

\noindent{\textbf{Definition 2.1 :}} A sequence $\{x_{n}\}_{n\in \mathbb{N}}$  of real numbers is said to be statistically convergent of order $\alpha$ (where $0<\alpha\leq 1$) to a real number $x$  if for each $\varepsilon>0,$ the set $K=\{n\in \mathbb{N}: |x_n-x|\geq \varepsilon\}$ has $\alpha$-natural density zero, i.e $$\displaystyle{\lim_{n\rightarrow \infty}}\frac{1}{n^{\alpha}}|\{k\leq n:|x_k-x|\geq \varepsilon\}|=0$$ and we write $x_{n}\xrightarrow{S^{\alpha}} x.$\\

Now we would like to introduce the definition of statistical convergence of order $\alpha$ in probability for a sequence of random variables as follows.\\

\noindent{\textbf{Definition 2.2.}} Let $(S,\triangle,P)$ be a probability space and $\{X_{n}\}_{n\in \mathbb{N}}$ be a sequence of random variables where each $X_{n}$ is defined on the same sample space $S$ (for each $n$) with respect to a given class of events $\triangle$ and a given probability function  $P:\triangle \rightarrow \mathbb{R}.$ Then the sequence $\{X_{n}\}_{n\in \mathbb{N}}$ is said to be statistically convergent of order $\alpha$ (where $0<\alpha\leq 1$) in probability to a random variable $X$ (where $X:S\rightarrow \mathbb{R}$) if for any $\varepsilon,\delta >0$ $$ {\lim_{n\rightarrow \infty}}\frac{1}{n^{\alpha}}|\{k\leq n: P(|X_{k}-X|\geq \varepsilon)\geq \delta\}|=0$$ or equivalantly $$ {\lim_{n\rightarrow \infty}}\frac{1}{n^{\alpha}}|\{k\leq n: 1-P(|X_{k}-X|< \varepsilon)\geq \delta\}|=0.$$ In this case we write $S^{\alpha}-\lim P(|X_{k}-X|\geq \varepsilon)=0$  or  $S^{\alpha}-\lim P(|X_{k}-X|< \varepsilon)=1$ or by $X_{n}\xrightarrow{PS^{\alpha}} X.$  The class of all  sequences of random variables which are statistically convergent of order $\alpha$ in probability is denoted simply by $PS^{\alpha}.$  \\

\noindent{\textbf{Theorem 2.1.}} If a sequence of constants $x_{n}\xrightarrow{S^{\alpha}} x$ then regarding a constant as a random variable having one point distribution at that point, we may also write $x_{n}\xrightarrow{PS^{\alpha}}x$.\\

\noindent{\textbf{Proof.}}  Let $\varepsilon >0$ be any arbitrarily small positive real number. Then $$\displaystyle{\lim_{n\rightarrow \infty}}~\frac{1}{n^{\alpha}}|\{k\leq n:|x_{k}-x|\geq \varepsilon\}|=0.$$ Now let $\delta>0$. So the set $\{k \in  \mathbb{N}:P(|x_{k}-x|\geq \varepsilon)\geq \delta\}\subseteq K$ where $K=\{k\in \mathbb{N}:|x_{k}-x|\geq \varepsilon\}$. This shows that  $x_{n}\xrightarrow{PS^{\alpha}} x$.\\

The following example shows that in general the converse of Theorem 2.1 is not true and also shows that there is a sequence $\{X_{n}\}_{n\in \mathbb{N}}$ of random variables which is statistically convergent in probability to a random variable X but it is not statistically convergent of order $\alpha$ in probability for $0<\alpha<1.$\\

\noindent{\textbf{Example 2.1.}} Let $\frac{r}{s}$ be a rational number between $\alpha$ and $\beta$. Let the probability density function of $X_{n}$ be given by,

If $n=[m^{\frac{s}{r}}]$ let, \begin{equation*}
~f_{n}(x)=
\begin{cases}
1  ~\mbox{where} ~ 0<x<1 ~&\\
0  ~\mbox{otherwise},
\end{cases}
\end{equation*}

If $n\neq [m^{\frac{s}{r}}]$ define,\begin{equation*}
f_{n}(x)=
\begin{cases}
\frac{nx^{n-1}}{2^{n}} ~ \mbox{where} ~~ 0<x<2~&\\
0 ~ \mbox{otherwise.}
\end{cases}
\end{equation*}
Now let $0<\varepsilon,\delta<1.$

Then \begin{equation*}
P(|X_{n}-2|\geq \varepsilon)=
\begin{cases}
1 ~ \mbox{if} ~ n=[m^{\frac{s}{r}}] ~ \mbox{for some} ~~  m\in \mathbb{N}~&\\
(1- \frac{\varepsilon}{2})^{n} ~ \mbox{if} ~ n\neq[m^{\frac{s}{r}}] ~~ \mbox{for any} ~~  m\in \mathbb{N}.
\end{cases}
\end{equation*}

Consequently we have the inequality, $$ \lim_{n\rightarrow \infty}{\{\frac{n^{\frac{r}{s}}-1}{n^{\alpha}}\}} \leq  \lim_{n\rightarrow \infty} \frac{1}{n^{\alpha}}|{\{k\leq n:P(|X_{k}-2|\geq \varepsilon)\geq \delta \}}|$$
and $$\lim_{n\rightarrow \infty} \frac{1}{n^{\beta}}|{\{k\leq n:P(|X_{k}-2|\geq \varepsilon)\geq \delta \}}| \leq  \lim_{n\rightarrow \infty}\frac{1}{n^{\beta}}\{n^{\frac{r}{s}}+1\}+\lim_{n\rightarrow \infty}\frac{c}{n^{\beta}}$$
(where c is a fixed finite positive integer). This shows that $X_{n}$ is statistically convergent of order $\beta$ in probability to $2$ but is not statistically convergent of order $\alpha$ in probability to $2$ for any $\alpha<\beta$ and this is not the usual statistical convergence of order $\alpha$ of real numbers. So the converse of Theorem 2.1 is not true. Also by taking $\beta =1,$ we see that $X_{n}\xrightarrow {PS}2$ but $\{X_{n}\}_{n\in \mathbb{N}}$ is not statistically convergent of order $\alpha$ in probability to $2$ for $0<\alpha<1.$\\

\noindent{\textbf{Theorem 2.2.}} [Elementary properties].\\
(i) If $X_{n}\xrightarrow {PS^{\alpha}} X$ and $X_{n}\xrightarrow {PS^{\beta}} Y$ then $P\{X=Y\}=1$ for any $\alpha,\beta$ where $0<\alpha,\beta \leq 1$\\
(ii) If $X_{n}\xrightarrow {PS^{\alpha}} X$ and $Y_{n}\xrightarrow {PS^{\beta}} Y$ then $X_{n}+Y_{n}\xrightarrow {PS^{\max\{\alpha,\beta\}}} X+Y$ for any $\alpha,\beta$ where $0<\alpha,\beta \leq 1$\\
(iii) If $X_{n}\xrightarrow {PS^{\alpha}} X$ then $cX_{n}\xrightarrow {PS^{\beta}} cX$ where $0<\alpha<\beta\leq 1.$\\
(iv) Let $0<\alpha\leq\beta\leq 1.$ Then $PS^{\alpha}\subset PS^{\beta}$ and this inclusion is strict whenever $\alpha<\beta$.\\


\noindent{\textbf{Proof.}} (i) Without loss of generality we assume $\beta\leq \alpha.$ If possible let $P\{X=Y\}\neq 1.$ Then there exists two positive real numbers $\varepsilon$ and $\delta$ such that $P(|X-Y|\geq \varepsilon)=\delta>0.$ Then we have
$$\displaystyle\lim_{n\rightarrow \infty}\frac{n}{n^{\alpha}}\leq \lim_{n\rightarrow \infty}\frac{1}{n^{\alpha}}|{\{k\leq n:P(|X_{k}-X|\geq \frac{\varepsilon}{2})\geq \frac{\delta}{2}}\}|+ \lim_{n\rightarrow \infty}\frac{1}{n^{\beta}}|{\{k\leq n:P(|X_{k}-Y|\geq \frac{\varepsilon}{2})\geq \frac{\delta}{2}}\}|$$ which is impossible because the left hand limit is not $0$ whereas the right hand limit is $0.$ So $P\{X=Y\}=1.$  \\

Proofs of (ii), (iii) are straightforward and so are omitted.\\

(iv) The first part is obvious. The inclusion is proper as can be seen from Example 2.1.\\

\noindent{\textbf{Remark 2.1.}  In Theorem 2 \cite{2}, we get, if $\alpha<\beta$ then $m^\alpha_0\subset m^\beta_0 $ (i.e., statistical convergence of order $\alpha$ $\Rightarrow$ statistical convergence of order $\beta$) and this inclusion is strict for at least those $\alpha, \beta$ for which there is a $k\in \mathbb{N}$ such that $\alpha< \frac{1}{k}<\beta.$ But the above Theorem 2.2 (iv) and Example 2.1 shows that the inclusion is strict for any $\alpha<\beta$  (i.e., $\alpha, \beta$ may satisfy the inequality $\frac{1}{k+1}<\alpha<\beta<\frac{1}{k}$ for any $k \in \mathbb{N}$).\\

\noindent{\textbf{Definition 2.3.} Let $(S,\triangle,P)$ be a probability space and $\{X_{n}\}_{n\in \mathbb{N}}$ be a sequence of random variables where each $X_{n}$ is defined on the same sample space $S$ (for each $n$) with respect to a given class of events $\triangle$ and a given probability function  $P:\triangle \rightarrow \mathbb{R}.$   A sequence of random variables ${\{X_{n}\}}_{n\in \mathbb{N}}$ is said to be strong $p$-Ces$\grave{\mbox{a}}$ro summable of order $\alpha$ (where $0<\alpha \leq 1$ and $p>0$ is any fixed positive real number) in probability to a random variable X  if for any $\varepsilon>0$ $$\lim_{n\rightarrow \infty}\frac{1}{n^{\alpha}}\displaystyle{\sum_{k=1}^{n}}{\{P(|X_{k}-X|\geq \varepsilon)\}}^{p}=0.$$ In this case we write $X_{n}\xrightarrow {PW_{p}^{\alpha}} X.$ The class of all sequences of random variables which are strong $p$-Ces$\grave{\mbox{a}}$ro summable of order $\alpha$ in probability is denoted simply by $PW_{p}^{\alpha}.$\\


\noindent{\textbf{Theorem 2.3.} (i) Let $0<\alpha\leq \beta \leq 1.$ Then $PW_{p}^{\alpha}\subset PW_{p}^{\beta}.$ This inclusion is strict whenever $\alpha<\beta.$\\
(ii) Let $0<\alpha\leq 1$ and $0<p<q<\infty$. Then $PW_{q}^{\alpha}\subset PW_{p}^{\alpha}.$\\

\noindent{\textbf{Proof.}} (i) The first part of this theorem is straightforward and so is omitted. For the second part we will give an example to show that there is a sequence of random variables ${\{X_{n}\}}_{n\in \mathbb{N}}$ which is strong p-Cesaro summable of order $\beta$ in probability to a random variable X but is not strong p-Cesaro summable of order $\alpha$ in probability for any $\alpha<\beta.$\\

Let $\frac{r}{s}$ be a rational number between $\alpha$ and $\beta$. We consider a sequence of random variables :
\begin{equation*}
X_{n} \in
\begin{cases}
 \{-1,1\}~ \mbox{with  probability} ~ \frac{1}{2},~
\mbox{if} ~ n=[m^\frac{s}{r}] ~~ \mbox{for some} ~ m\in \mathbb{N}~&\\
 \{0,1\} ~~ \mbox{with probability} ~~ P(X_{n}=0)=1-\frac{1}{\sqrt[p]{n^2}} ~~ \mbox{and} ~~ P(X_{n}=1)=\frac{1}{\sqrt[p]{n^2}}, ~ \mbox{if} ~~ n\neq [m^\frac{s}{r}] ~~ \mbox{for any} ~ m\in \mathbb{N}
\end{cases}
\end{equation*}


For $0<\varepsilon<1$ we get the inequality, $$\lim_{n\rightarrow \infty}\frac{n^\frac{r}{s}-1}{n^{\alpha}}\leq \lim_{n\rightarrow \infty}\frac{1}{n^{\alpha}}\displaystyle{\sum_{k=1}^{n}}{\{P(|X_{k}-0|\geq \varepsilon)}\}^{p}$$ and
$$\lim_{n\rightarrow \infty}\frac{1}{n^{\beta}}\displaystyle{\sum_{k=1}^{n}}{\{P(|X_{k}-0|\geq \varepsilon)}\}^{p}\leq \lim_{n\rightarrow \infty} [\frac{n^\frac{r}{s}+1}{n^{\beta}} + \frac{1}{n^{\beta}}(\frac{1}{1^{2}}+\frac{1}{2^{2}}+...+\frac{1}{n^{2}})]$$
This shows that $X_{n}\xrightarrow {PW_{p}^{\beta}} 0$ but $\{X_{n}\}_{n\in \mathbb{N}}$ is not strong $p$-Cesaro summable of order $\alpha$ in probability to $0$.\\

(ii) Proof is straightforward and so is omitted.\\

\noindent{\textbf{Theorem 2.4.} Let $0<\alpha \leq \beta \leq 1.$ Then $PW_{p}^{\alpha}\subset PS^{\beta}.$\\

\noindent{\textbf{Proof.}}~ For any $\varepsilon,\delta>0,$
$$\displaystyle{\sum_{k=1}^{n}}{\{P(|X_{k}-X|\geq \varepsilon)}\}^{p}  \geq
\underset{P(|X_{k}-X|\geq \varepsilon)\geq \delta} {\underset{ k=1 }{\sum^{n}}}{\{P(|X_{k}-X|\geq \varepsilon)}\}^{p}
   \geq       \delta^{p} |{\{k \leq n :P(|X_{k}-X|\geq \varepsilon)\geq \delta}\}|$$

$\Rightarrow~{\underset{n\rightarrow \infty}{\lim}} \frac{1}{n^{\alpha}}\displaystyle{\sum_{k=1}^{n}}{\{P(|X_{k}-X|\geq \varepsilon)}\}^{p}\geq \displaystyle{\lim_{n\rightarrow \infty}}\frac{1}{n^{\beta}} |{\{k \leq n :P(|X_{k}-X|\geq \varepsilon)\geq \delta}\}|.\delta^{p}$\\
This proves the theorem.\\

\noindent{\textbf{Note 2.1.}  If a sequence of random variables ${\{X_{n}}\}_{n\in \mathbb{N}}$ is strong $p$-Ces$\grave{\mbox{a}}$ro summable of order $\alpha$ in probability to X then it is statistically convergent of order $\alpha$ in probability to X i.e $PW_{p}^{\alpha}\subset PS^{\alpha}.$\\

But the converse of Theorem 2.4 (or Note 2.1) is not generally true as can be seen from the following example.\\

\noindent{\textbf{Example 2.2.} Let a sequence of random variable ${\{X_{n}}\}_{n\in \mathbb{N}}$ be defined by,\\
\begin{equation*}
X_{n}\in
\begin{cases}
\{-1,1\} ~~\mbox{with probability} ~~ \frac{1}{2},~~ \mbox{if} ~~ n=m^{m} ~~ \mbox{for some}~~ m\in \mathbb{N}~&\\
\{0,1\} ~~ \mbox{with probability}~~ P(X_{n}=0)=1-\frac{1}{\sqrt[2p]{n}},~~ P(X_{n}=1)=\frac{1}{\sqrt[2p]{n}},~~ \mbox{if}~ n\neq m^{m} ~ \mbox{for any} ~ m\in \mathbb{N}
\end{cases}
\end{equation*}
Let $0<\varepsilon<1$ be given.

Then \begin{equation*}
P(|X_{n}-0|\geq \varepsilon)=
\begin{cases}
1 ~~  \mbox{if} ~~ n=m^{m} ~~\mbox{for some}~~ m\in \mathbb{N}~&\\
\frac{1}{\sqrt[2p]{n}} ~~  \mbox{if} ~~ n\neq m^{m} ~~\mbox{for any}~~ m\in \mathbb{N}.
\end{cases}
\end{equation*}

This implies $X_{n}\xrightarrow {PS^{\alpha}} 0$ for each $0<\alpha\leq 1.$ Now let $H={\{n\in \mathbb{N}: ~~ n\neq m^{m} ~~ \mbox{for any}~ m\in \mathbb{N}}\}.$
$$\mbox{Then}~ \frac{1}{n^{\alpha}}\displaystyle{\sum_{k=1}^{n}}{\{P(|X_{k}-0|\geq \varepsilon)}\}^{p}=\frac{1}{n^{\alpha}}\underset{k\in H}{\underset{k=1}{\sum^{n}}}{\{P(|X_{k}-0|\geq \varepsilon)}\}^{p}+\frac{1}{n^{\alpha}}\underset{k\notin H}{\underset{k=1}{\sum^{n}}}{\{P(|X_{k}-0|\geq \varepsilon)}\}^{p}$$
$$= \frac{1}{n^{\alpha}}\underset{k\in H}{\underset{k=1}{\sum^{n}}}\frac{1}{\sqrt{k}} +\frac{1}{n^{\alpha}}\underset{k\notin H}{\underset{k=1}{\sum^{n}}}1 >\frac{1}{n^{\alpha}}{\underset{k=1}{\sum^{n}}}\frac{1}{\sqrt{k}}> \frac{1}{n^{\alpha - \frac{1}{2}}} ~~(\mbox{since}~~{\underset{k=1}{\sum^{n}}}\frac{1}{\sqrt{k}}>\sqrt{n}~~ \mbox{for}~~ n\geq 2).$$
So $X_{n}$ is not strong $p$-Ces$\grave{\mbox{a}}$ro summable of order $\alpha$ in probability  to $0$ for $0<\alpha\leq \frac{1}{2}$.\\

\noindent{\textbf{Theorem 2.5.} For $\alpha=1,$ $PW_{p}^{\alpha}=PS^{\alpha}$ i.e in other words $PW_{p}=PS$\\

\noindent{\textbf{Proof.}} That $PW_{p}\subset PS$ readily follows from Theorem 2.4. For the converse part take $\varepsilon,\delta>0.$ Now
$$\frac{1}{n}\displaystyle{\sum_{k=1}^{n}}{\{P(|X_{k}-X|\geq \varepsilon)}\}^{p}=\frac{1}{n}\underset{P(|X_{k}-X|\geq \epsilon)\geq \delta}{\underset{k=1}{\sum^{n}}}{\{P(|X_{k}-X|\geq \varepsilon)}\}^{p}+\frac{1}{n}\underset{P(|X_{k}-X|\geq \varepsilon)< \delta}{\underset{k=1}{\sum^{n}}}{\{P(|X_{k}-X|\geq \varepsilon)}\}^{p}$$
$$\leq \frac{1}{n}|{\{k\leq n:P(|X_{k}-X|\geq \varepsilon)\geq \delta}\}|+ \frac{1}{n}.n.\delta^{p}.$$
This shows that $PS\subset PW_{p}$. Hence the result follows.\\

\section{\textbf{Lacunary Statistical Convergence of Order $\alpha$ in Probability}}

First of all we would like to recall the definition of lacunary statistical convergence of order $\alpha$ for a sequence of real numbers from \cite{E1} as follows:\\

\noindent{\textbf{Definition 3.1.} A real sequence $\{x_{n}\}_{n\in \mathbb{N}}$ is said to be lacunary statistically convergent of order $\alpha$ to a real number $L$ if for any $\varepsilon>0$ $$\displaystyle{\lim_{r\rightarrow \infty}} \frac{1}{h_{r}^{\alpha}}|\{k\in I_{r}:|x_{k}-L|\geq \varepsilon\}|=0.$$  where $\theta=\{k_{r}\}_{r\in \mathbb{N}\cup \{0\}}$ be the lacunary sequence (for definition of lacunary sequence see \cite{9}), $h_{r}=(k_{r}-k_{r-1}),$  $I_{r}=(k_{r-1},k_{r}]$ and $q_r=\frac{k_r}{k_{r-1}}.$ In this case we write $x_{n}\xrightarrow {S_{\theta}^{\alpha}}L$.\\

 Now we would like to introduce the definition of lacunary statistical convergence of order $\alpha$ in probability for a sequence of random variables as follows.\\

\noindent{\textbf{Definition 3.2.} Let $\theta=\{k_{r}\}_{r\in \mathbb{N}\cup \{0\}}$ be a lacunary sequence and $\{X_{n}\}_{n\in \mathbb{N}}$ be a sequence of random variables where each $X_{n}$ is defined on the same sample space $S$ (for each $n$) with respect to a given class of events $\triangle$ and a given probability function  $P:\triangle \rightarrow \mathbb{R}$. Then the sequence ${\{X_{n}}\}_{n\in \mathbb{N}}$ is said to be lacunary statistically convergent or $S_{\theta}$-convergent of order $\alpha$(where $0<\alpha\leq 1$) in probability to a random variable X if for any $\varepsilon,\delta>0,$ $$\displaystyle{\lim_{r\rightarrow \infty}} \frac{1}{h_{r}^{\alpha}}|{\{k\in I_{r}:P(|X_{k}-X|\geq \varepsilon)\geq \delta}\}|=0.$$ In this case we write $X_{n}\xrightarrow{PS_{\theta}^{\alpha}} X.$ The class of all sequences of random variables which are $S_{\theta}$-convergent of order $\alpha$ in probability is simply denoted by $PS_{\theta}^{\alpha}.$\\

\noindent{\textbf{Definition 3.3.} Let $\theta=\{k_{r}\}_{r\in \mathbb{N}\cup \{0\}}$ be a lacunary sequence and $\{X_{n}\}_{n\in \mathbb{N}}$ be a sequence of random variables where each $X_{n}$ is defined on the same sample space $S$ (for each $n$) with respect to a given class of events $\triangle$ and a given probability function  $P:\triangle \rightarrow \mathbb{R}$. The sequence ${\{X_{n}}\}_{n\in \mathbb{N}}$ is said to be $N_{\theta}$-convergent of order $\alpha$ (where $0<\alpha\leq 1$) in probability  to a random variable X if for any $\varepsilon>0,$  $$\displaystyle{\lim_{r\rightarrow \infty}}\frac{1}{h_{r}^{\alpha}}\sum_{k\in I_{r}}P(|X_{k}-X|\geq \varepsilon)=0.$$
In this case we write $X_{n}\xrightarrow{PN_{\theta}^{\alpha}} X.$ The class of all sequences of random variables which are $N_{\theta}$-convergent of order $\alpha$ in probability is simply denoted by $PN_{\theta}^{\alpha}.$\\

\noindent{\textbf{Theorem 3.1.} For a fixed $\theta,$ $PS_{\theta}^{\alpha}$-limit and  $PN_{\theta}^{\alpha}$-limit of a sequence of random variables ${\{X_{n}}\}_{n\in \mathbb{N}}$ are unique.\\

The proof is omitted.\\

However the $PS_{\theta}$-limit of a sequence of random variables for two different lacunary sequences may not be equal as can be seen from the following example.\\

\noindent{\textbf{Example 3.1.} Let $G:\mathbb{R}^+\rightarrow \mathbb{N}$ be defined by,
$$G(x)= n~~ \mbox{if}~~ n!<x\leq (n+1)! $$
Let us take two lacunary sequences $\theta_1=\{(2r)!\}$ ~ and ~ $\theta_2=\{(2r+1)!\}$. Let us define a sequence of random variables $\{X_{n}\}_{n\in \mathbb{N}}$ by,
\begin{equation*}
X_{n}\in
\begin{cases}
\{-1,1\}~~ \mbox{with probability}~~ P(X_{n}=-1)=\frac{1}{n}, P(X_{n}=1)=1-\frac{1}{n},~~ \mbox{if}~~G(n)~ \mbox{is even}, ~&\\
 \{0,1\}~~ \mbox{with probability}~~ P(X_{n}=0)=1- \frac{1}{n},P(X_{n}=1)=\frac{1}{n},~~ \mbox{if}~~G(n)~ \mbox{is odd},
\end{cases}
\end{equation*}
Let  $0<\varepsilon, \delta< 1.$ For the lacunary sequence $\theta_{1},$
$$\frac{1}{h_{r+1}}|\{k\in I_{r+1}:P(|X_{k}-0|\geq \varepsilon)\geq \delta\}|=\frac{1}{h_{r+1}}[\{(2r+1)!-(2r)!\}+c_1]\rightarrow 0 ~~ as ~~ r\rightarrow \infty$$ where $c_1$ is a finite positive integer. This shows that $X_{n}\xrightarrow {PS_{\theta_1}} 0.$\\

For the lacunary sequence $\theta_2,$
$$\frac{1}{h_{r+1}}|\{k\in I_{r+1}:P(|X_{k}-1|\geq \epsilon)\geq \delta\}|\leq \frac{\{(2r+2)!-(2r+1)!\}+c_2}{(2r+3)!-(2r+1)!}\rightarrow 0 ~~ as~~ r\rightarrow \infty$$ where $c_2$ is a finite positive integer. This shows that $X_{n}\xrightarrow {PS_{\theta_2}} 1.$\\

\noindent{\textbf{Remark 3.1.}} As in Theorem 2.5, we can easily establish that $PN_{\theta}=PS_{\theta}$.\\

\noindent{\textbf{Theorem 3.2.}} Let $0<\alpha\leq \beta \leq 1.$ Then $PN_{\theta}^{\alpha}\subset PS_{\theta}^{\beta}.$ This inclusion is strict whenever $\alpha<\beta.$\\

\noindent{\textbf{Proof.}} Let $\varepsilon,\delta>0.$ Now $$\displaystyle{\sum_{k\in I_{r}}}P(|X_{k}-X|\geq \varepsilon)\geq \underset{P(|X_{k}-X|\geq \varepsilon)\geq \delta}{\underset{k\in I_{r}}\sum}P(|X_{k}-X|\geq \varepsilon)\geq \delta |{\{k\in I_{r}:P(|X_{k}-X|\geq \varepsilon)\geq \delta}\}|$$
which implies that,
\begin{eqnarray*}
\frac{1}{h_{r}^{\alpha}}\displaystyle{\sum_{k\in I_{r}}}P(|X_{k}-X|\geq \varepsilon)
&\geq& \frac{\delta}{h_{r}^{\alpha}}|{\{k\in I_{r}:P(|X_{k}-X|\geq \varepsilon)\geq \delta}\}|\\
&\geq& \frac{\delta}{h_{r}^{\beta}}|{\{k\in I_{r}:P(|X_{k}-X|\geq \varepsilon)\geq \delta}\}|.
\end{eqnarray*}
As $r\rightarrow \infty$ we have, $ \frac{1}{h_{r}^{\beta}}|{\{k\in I_{r}:P(|X_{k}-X|\geq \varepsilon)\geq \delta}\}|\rightarrow 0$. Hence the result follows.\\
For the second part of this theorem we will constract an example.\\

Let $c$ be a rational number between $\alpha$ and $\beta.$ Let a sequence of random variables ${\{X_{n}}\}_{n\in \mathbb{N}}$ be defined by
\begin{equation*}
X_{n}\in
\begin{cases}
\{-1,1\} ~~ \mbox{with probability} ~~ \frac{1}{2},~~ \mbox{if} ~~ n ~\mbox{is the first}~~[h_{r}^{c}]~\mbox{integers in the interval}~I_{r}, ~&\\
 \{0,1\} ~~ \mbox{with probability} ~ P(X_{n}=0)=(1-\frac{1}{n})~~ \mbox{and} ~~ P(X_{n}=1)=\frac{1}{n},~~\mbox{if}~~ n ~\mbox{is other than the first}\\
~~[h_{r}^{c}]~~ \mbox{integers in the interval}~~ I_{r}.
\end{cases}
\end{equation*}

Now let $0<\varepsilon<1$ and $\delta<1.$ Then we have,
\begin{equation*}
P(|X_{n}-0|\geq \varepsilon)=
\begin{cases}
1 ~~ \mbox{if} ~~n ~\mbox{is the first}~~[h_{r}^{c}]~\mbox{integers in the interval}~I_{r}, ~&\\
\frac{1}{n}~~\mbox{if}~~ n ~\mbox{is other than the first}~~[h_{r}^{c}]~\mbox{integers in the interval}~I_{r}.
\end{cases}
\end{equation*}

Now $\frac{1}{h_{r}^{\beta}}|{\{k\in I_{r}:P(|X_{k}-0|\geq \varepsilon)\geq \delta}\}|\leq \frac{1}{h_{r}^{\beta}}\{h_{r}^{c}+d\}=(\frac{1}{h_{r}^{\beta -c}}+\frac{d}{h_{r}^{\beta}})\rightarrow 0 ~~\mbox{as}~~ r\rightarrow \infty$ since $c<\beta\leq1$\\
where $d$ is a finite positive integer. So $X_{n}\xrightarrow {PS_{\theta}^{\beta}} 0$.\\ We also have the inequality
$$\frac{h_{r}^{c}-1}{h_{r}^{\alpha}}\leq \frac{1}{h_{r}^{\alpha}}\sum_{k\in I_{r}}P(|X_{k}-0|\geq \varepsilon)$$\\
which shows that $X_{n}$ is not $N_{\theta}$-convergent of order $\alpha$ in probability to $0$ since $0<\alpha<c$. So the inequality is strict whenever $\alpha<\beta$.\\

\noindent{\textbf{Note 3.1.}} We show that for bounded sequences the equality between the two sets $PN_{\theta}^{\alpha}$ and $PS_{\theta}^{\beta}$ may not be true in general whenever $\alpha<\beta.$\\

\noindent{\textbf{Theorem 3.3.}} If $0<\alpha\leq \beta\leq 1$ then $PS^{\alpha}\subset PS_{\theta}^{\beta}$ if and only if $\liminf q_{r}>1.$ If $\liminf q_{r}=1$ then there exists a sequence of random variables $\{X_{n}\}_{n\in \mathbb{N}}$ which is not $S_{\theta}$-convergent of order $\beta$ in probability but is statistically convergent of order $\alpha$ in probability for some $\alpha,\beta.$\\

\noindent{\textbf{Proof.}} Let $\liminf q_{r}>1 ~\mbox{and}~ X_{n}\xrightarrow {PS^{\alpha}}X.$ Then for each $\delta>0$ we can find $q_{r}$ for sufficiently large $r$ such that $q_{r}\geq (1+\delta)\Rightarrow \frac{h_{r}}{k_{r}}\geq \frac{\delta}{(1+\delta)}\Rightarrow \frac{1}{h_{r}^{\alpha}}\leq \frac{1}{k_{r}^{\alpha}} \frac{(1+\delta)^{\alpha}}{\delta^{\alpha}}.$
Let $\varepsilon>0.$
\begin{eqnarray*}
\frac{1}{h_{r}^{\beta}}|\{k\in I_{r}:P(|X_{k}-X|\geq \varepsilon)\geq \delta)\}|
&\leq&\frac{1}{h_{r}^{\alpha}}|\{k\leq k_{r}:P(|X_{k}-X|\geq \varepsilon)\geq \delta)\}|\\
&\leq& \frac{(1+\delta)^{\alpha}}{\delta^{\alpha}}\frac{1}{k_{r}^{\alpha}}|\{k\leq k_{r}:P(|X_{k}-X|\geq \varepsilon)\geq \delta)\}|.
\end{eqnarray*}
This shows that $X_{n}\xrightarrow {PS_{\theta}^{\beta}}X.$\\

Now let $PS^{\alpha}\subset PS_{\theta}^{\beta}(\mbox{where}~ 0<\alpha\leq \beta\leq 1).$ If possible let $\liminf q_{r}=1.$ From \cite{F2} we can find a subsequence $k_{r(j)}$ satisfying $\frac{k_{r(j)}}{k_{r(j)-1}}<(1+\frac{1}{j}) ~\mbox{and}~ \frac{k_{r(j)-1}}{k_{r(j-1)}}>j~\mbox{where}~j\in \mathbb{N}$\\

We define a sequence of random variables by
\begin{equation*}
X_{n}\in
\begin{cases}
 \{-1,1\}~\mbox{with probability}~ \frac{1}{2}~~\mbox{if}~ n\in I_{r(j)}~~\mbox{where}~ j\in \mathbb{N}~&\\
 \{0,1\}~\mbox{with probability}~ P(X_{n}=0)=1-\frac{1}{n^{2}}~\mbox{and}~P(X_{n}=1)=\frac{1}{n^{2}}~~\mbox{if}~ n\notin I_{r(j)}~~\mbox{for any}~ j\in \mathbb{N}
\end{cases}
\end{equation*}
Let $0<\varepsilon<1~ \mbox{and}~ \delta<1.$ Now
$$P(|X_{n}-0|\geq \varepsilon)=1 ~~\mbox{if}~ n\in I_{r(j)}~~\mbox{where}~ j\in \mathbb{N},$$ and
$$P(|X_{n}-0|\geq \varepsilon)=\frac{1}{n^{2}} ~~\mbox{if}~ n\notin I_{r(j)}~~\mbox{for any} j\in \mathbb{N}.$$\\
Now $\frac{1}{h_{r(j)}^{\beta}}|\{k\in I_{r(j)}:P(|X_{n}-0|\geq \varepsilon)\geq \delta\}|=\frac{1}{h_{r(j)}^{\beta}}h_{r(j)}\rightarrow \infty~ \mbox{as}~r\rightarrow \infty.$ As $\frac{1}{h_{r(j)}^{\beta}}|\{k\in I_{r(j)}:P(|X_{n}-0|\geq \varepsilon)\geq \delta\}|$ is a subsequence of the sequence $\frac{1}{h_{r}^{\beta}}|\{k\in I_{r}:P(|X_{n}-0|\geq \varepsilon)\geq \delta\}|,$ this shows that $X_{n}$ is not $S_{\theta}$-convergent of order $\beta$(where $0<\beta\leq 1$)in probability to $0$.\\

Finally let $\alpha=1.$  If we take t sufficiently large such that $k_{r(j)-1}<t\leq k_{r(j+1)-1}$ then we observe that,\\
$\frac{1}{t}\displaystyle{\sum_{k=1}^{t}}P(|X_{k}-0|\geq \varepsilon)\leq \frac{k_{r(j-1)}+h_{r(j)}}{k_{(r(j)-1)}}+\frac{1}{t}\{1+ \frac{1}{2^{2}}+...+\frac{1}{t^{2}}\}\leq \frac{2}{j}+\frac{1}{t}\{1+ \frac{1}{2^{2}}+...+\frac{1}{t^{2}}\}\rightarrow 0~ \mbox{if}~ j,t\rightarrow \infty$\\
This shows that $X_{n}\xrightarrow {PS}0$ and we conclude that $\liminf q_{r}~\mbox{must be}>1.$\\

\noindent{\textbf{Theorem 3.4.}} Let $X_{n}\xrightarrow {PS^{\alpha}} X$ and $X_{n}\xrightarrow {PS_{\theta}^{\beta}}Y$. If $\liminf q_{r}>1$ and $0<\beta\leq \alpha\leq 1$ then $P\{X=Y\}=1.$\\

\noindent{\textbf{Proof.}} Let $\varepsilon> 0$ be any small positive real number and if possible let $P(|X-Y|\geq \varepsilon)=\delta>0.$ ~~ Now we have the inequality $\{P(|X-Y|\geq \varepsilon)\}\leq \{P(|X_{n}-X|\geq \frac{\varepsilon}{2})\}+ \{P(|X_{n}-Y|\geq \frac{\varepsilon}{2})\}.$\\

So ${\{k\in I_{r}:P(|X-Y|\geq \varepsilon)\geq \delta}\} \subseteq {\{k\in I_{r}: P(|X_{k}-X|\geq \frac{\varepsilon}{2})\geq \frac{\delta}{2}}\} \bigcup {\{k\in I_{r}: P(|X_{k}-Y|\geq \frac{\varepsilon}{2})\geq \frac{\delta}{2}}\}$\\

$\Rightarrow |{\{k\in I_{r}:P(|X-Y|\geq \varepsilon)\geq \delta}\}|\leq |{\{k\in I_{r}: P(|X_{k}-X|\geq \frac{\varepsilon}{2})\geq \frac{\delta}{2}}\}|+ |{\{k\in I_{r}: P(|X_{k}-Y|\geq \frac{\varepsilon}{2})\geq \frac{\delta}{2}}\}|$\\

$\Rightarrow |{\{k\in I_{r}:P(|X-Y|\geq \varepsilon)\geq \delta}\}|\leq |{\{k\leq k_{r}: P(|X_{k}-X|\geq \frac{\varepsilon}{2})\geq \frac{\delta}{2}}\}|+ |{\{k\in I_{r}: P(|X_{k}-Y|\geq \frac{\varepsilon}{2})\geq \frac{\delta}{2}}\}|$\\

$\Rightarrow h_{r}\leq |{\{k\leq k_{r}: P(|X_{k}-X|\geq \frac{\varepsilon}{2})\geq \frac{\delta}{2}}\}|+ |{\{k\in I_{r}: P(|X_{k}-Y|\geq \frac{\varepsilon}{2})\geq \frac{\delta}{2}}\}|$\\

$\Rightarrow \{\frac{h_{r}}{k_{r}}\}^{\alpha}\leq \frac{1}{k_{r}^{\alpha}}|{\{k\leq k_{r}: P(|X_{k}-X|\geq \frac{\varepsilon}{2})\geq \frac{\delta}{2}}\}|+ \frac{1}{k_{r}^{\beta}}|{\{k\in I_{r}: P(|X_{k}-Y|\geq \frac{\varepsilon}{2})\geq \frac{\delta}{2}}\}|$\\

$\Rightarrow \{\frac{h_{r}}{k_{r}}\}^{\alpha}\leq \frac{1}{k_{r}^{\alpha}}|{\{k\leq k_{r}: P(|X_{k}-X|\geq \frac{\varepsilon}{2})\geq \frac{\delta}{2}}\}|+ \{\frac{h_{r}}{k_{r}}\}^{\beta}\frac{1}{h_{r}^{\beta}}|{\{k\in I_{r}: P(|X_{k}-Y|\geq \frac{\varepsilon}{2})\geq \frac{\delta}{2}}\}|$\\

$\Rightarrow \{1-\frac{1}{q_{r}}\}^{\alpha}\leq \frac{1}{k_{r}^{\alpha}}|{\{k\leq k_{r}: P(|X_{k}-X|\geq \frac{\varepsilon}{2})\geq \frac{\delta}{2}}\}|+ \{1-\frac{1}{q_{r}}\}^{\beta}\frac{1}{h_{r}^{\beta}}|{\{k\in I_{r}: P(|X_{k}-Y|\geq \frac{\varepsilon}{2})\geq \frac{\delta}{2}}\}|$\\

Taking $r\rightarrow \infty$ on both sides we see that the left hand side does not tend to zero since $\liminf q_{r}>1$ but the right hand side tends to zero. This is a contradiction. So we must have $P\{X=Y\}=1.$\\

\noindent{\textbf{Acknowledgments:}} We Thank the respected Referees for his/her careful reading of the paper and for several valuable suggestions which have improved the presentation of this paper.\\



\end{document}